\documentclass{article}
\usepackage{amssymb,amsmath,graphicx,ucs,hyperref}
\usepackage[utf8x]{inputenc}

\newtheorem{proposition}{Proposition}

\newtheorem{remark}{Remark}

\newtheorem{corollary}{Corollary}

\newenvironment{proof}{\noindent \emph{Proof. }}{\hfill \hbox{\rlap{$\sqcap$}$\sqcup$}\\}

\title{Density of Binary Disc Packings:\\
Lower and Upper Bounds}

\author{Thomas Fernique\footnote{LIPN, Univ. Paris Nord \& CNRS, fernique@lipn.fr}}
\date{}

\begin{document}
\maketitle

\begin{abstract}
We provide, for any $r\in (0,1)$, lower and upper bounds on the maximal density of a packing in the Euclidean plane of discs of radius $1$ and $r$.
The lower bounds are mostly folk, but the upper bounds improve the best previously known ones for any $r\in[0.11,0.74]$.
For many values of $r$, this gives a fairly good idea of the exact maximum density.
In particular, we get new intervals for $r$ which does not allow any packing more dense that the hexagonal packing of equal discs.
\end{abstract}

\noindent The computer programs referred to in this article can be found at the url:
\begin{center}
\url{https://arxiv.org/src/2107.14079v2/anc}
\end{center}
\section{Introduction}

A {\em disc packing} (or {\em circle packing}) is a set of interior-disjoint discs in the Euclidean plane.
Its {\em density} $\delta$ is the proportion of the plane covered by the discs:
$$
\delta:=\limsup_{k\to \infty}\frac{\textrm{area of the square $[-k,k]^2$ covered by discs}}{\textrm{area of the square $[-k,k]^2$}}.
$$

If all the discs have the same radius, it has been proven by T\'oth \cite{FT43} (see also \cite{CW10}) that the density is at most
$$
\delta_1:=\frac{\pi}{2\sqrt{3}}\approx 0.9069,
$$
reached for the so-called {\em hexagonal compact packing}, where discs are centered on a triangular grid (of size twice the disc radius).

What about more sizes of discs?
In particular, what about the maximal density of {\em binary disc packings}, that is, packings with discs of two different sizes?
Up to scaling, one can always assume that the largest disc has radius $1$ and the smallest one has radius $r\in (0,1)$.
We then denote by $\delta(r)$ the maximal density of packings by discs of radius $1$ and $r$.
Can we find the exact value of $\delta(r)$ for each $r\in(0,1)$?

The first decisive step was done by Blind \cite{Bli69}, who proved that $\delta(r)=\delta_1$ for $r\geq r_B$ where
$$
r_B:=\sqrt{\frac{7\tan\tfrac{\pi}{7}-6\tan\tfrac{\pi}{6}}{6\tan\tfrac{\pi}{6}-5\tan\tfrac{\pi}{5}}}\approx 0.74299.
$$
In other words, binary disc packings cannot achieve higher densities than packings of equal discs when the disc sizes are too close.
This yields the exact maximal density over the whole interval $[r_B,1)$.
Besides this, to the best of our knowledge, there are only $9$ specific ratios for which the exact value of $\delta(r)$ is known.
They are the ratios which allow {\em triangulated binary packings}, that is, packings with two sizes of discs whose contact graphs are triangulated (the vertices of the contact graph of a packing are the disc center, and the edges connect the centers of tangent discs).
These ratios are algebraic numbers which have been characterized in \cite{Ken06}.
The maximal density for each of them is proven in \cite{BF22} (see also \cite{Hep00,Hep03,Ken04} for the first cases).

If we cannot obtain an exact value for $\delta(r)$, can we find lower and upper bounds?
In particular, for which $r$ do we have $\delta(r)>\delta_1$, that is, which ratios allow binary disc packings which achieve higher densities than packings of equal discs?
This is a question that is of particular interest in materials science because the maximization of density seems to be a criterion for the formation of materials (due to attractive forces of the van der Waals type).
A ratio $r$ such that $\delta(r)>\delta_1$ thus suggests the possibility to obtain new materials by combining atoms or nanoparticles within this ratio.
A striking example is given by the $4$ binary assemblies of nanoparticles whose synthesis is explained in \cite{PDKM15}: they all correspond faithfully to one of the $9$ triangulated binary packings that have been proven to maximize the density.

\begin{figure}[hbt]
\centering
\includegraphics[width=\textwidth]{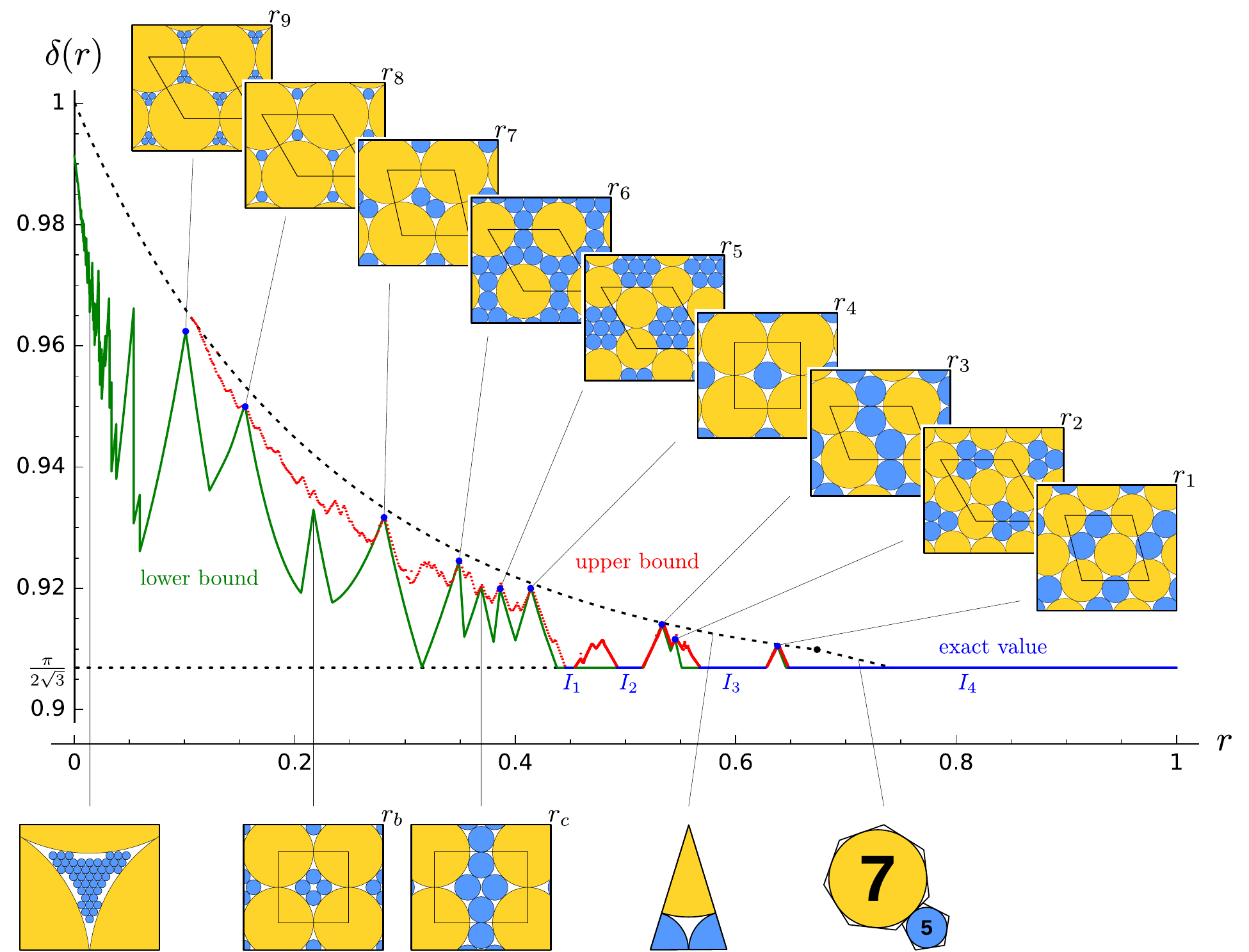}
\caption{
Lower bound (green curve) and upper bound (red curve) on $\delta(r)$.
The blue dots and lines indicate where lower and upper bounds coincide (the dots correspond to the $9$ triangulated packings, the lines to the intervals $I_i$'s in Cor.~\ref{cor:upper_bound}).
The best previous upper bound is indicated by the black dotted curve.
Labels $r_1$ through $r_9$, $r_b$ and $r_c$ denote disc ratios which appear in the text and are given in Table~\ref{tab:ratios}.
Zooms can be found in Fig.~\ref{fig:J34}, \ref{fig:J23}, \ref{fig:J12} and \ref{fig:J0}.
}
\label{fig:density}
\end{figure}

To obtain a lower bound, it suffices to find a ``good'' packing.
For example, if $r$ is small enough, namely $r\leq \tfrac{2}{\sqrt{3}}-1\approx 0.1547$, we can simply insert a small disc in each hole of a hexagonal compact packing of large discs to get $\delta(r)>\delta_1$.
More interesting, we can get lower bounds over a whole interval of ratios by continuously modifying a ``good'' packing so that the density decreases as little as possible.
A particularly efficient way to proceed is the so-called ``flipping and flowing method'', first used by T\'oth \cite{FT64}, see also \cite{CP19,CG20}.
In Section~\ref{sec:lower} we detail lower bounds obtained in this way.
These lower bounds seem to be mostly ``folk'' (see, {\em e.g.}, \cite{FJFS20,KenWeb}), but we provide here a SageMath worksheet \cite{sage} (the file \verb+lower_bounds.sage+ in supplementary materials) which give explicitly the transformations as well as the corresponding densities.
They are depicted in Fig.~\ref{fig:density} (green curve).
As a corollary, we get:

\begin{corollary}
\label{cor:lower_bound}
One has $\delta(r)>\delta_1$ for any $r$ in $(0,a_1)\cup (a_2,a_3)\cup (a_4,a_5)$ where $a_1\approx 0.4378$, $a_2\approx 0.5165$, $a_3\approx 0.5510$, $a_4\approx 0.6276$ and $a_5\approx 0.6456$ are algebraic number whose minimal polynomials are in the file \verb+lower_bounds.sage+.
\end{corollary}

\begin{table}
\centering
\begin{tabular}{|c|c|r|}
\hline
ratio & approx. value & minimal polynomial\\
\hline
$r_1$ & $0.6375559772$ & $x^4 \textrm{-} 10x^2 \textrm{-} 8x \textrm{+} 9$\\
$r_a$ & $0.6199144044$ & $x^4 \textrm{-} 12x^3 \textrm{-} 2x^2 \textrm{+} 4x \textrm{+} 1$\\
$r_2$ & $0.5451510421$ & $x^8 \textrm{-} 8x^7 \textrm{-} 44x^6 \textrm{-} 232x^5 \textrm{-} 482x^4 \textrm{-} 24x^3 \textrm{+} 388x^2 \textrm{-} 120x \textrm{+} 9$\\
$r_3$ & $0.5332964167$ & $8x^3 \textrm{+} 3x^2 \textrm{-} 2x \textrm{-} 1$\\
$r_4$ & $0.4142135624$ & $x^2 \textrm{+} 2x \textrm{-} 1$\\
$r_5$ & $0.3861061049$ & $9x^4 \textrm{-} 12x^3 \textrm{-} 26x^2 \textrm{-} 12x \textrm{+} 9$\\
$r_b$ & $0.3691023862$ & $x^3 \textrm{-} 5x^2 \textrm{-} x \textrm{+} 1$\\
$r_6$ & $0.3491981862$ & $x^4 \textrm{-} 28x^3 \textrm{-} 10x^2 \textrm{+} 4x \textrm{+} 1$\\
$r_7$ & $0.2807764064$ & $2x^2 \textrm{+} 3x \textrm{-} 1$\\
$r_c$ & $0.2168453354$ & $x^4 \textrm{-} 4x^3 \textrm{-} 2x^2 \textrm{-} 4x \textrm{+} 1$\\
$r_8$ & $0.1547005384$ & $3x^2 \textrm{+} 6x \textrm{-} 1$\\
$r_9$ & $0.1010205144$ & $x^2 \textrm{-} 10x \textrm{+} 1$\\
\hline
\end{tabular}
\label{tab:ratios}
\caption{Values of disc ratios appearing in Fig.~\ref{fig:density}, \ref{fig:flow841}, \ref{fig:flow873} and \ref{fig:flow9cba}.}
\end{table}

Finding upper bounds is more challenging.
For a given $r$, it indeed amouts to proving that among the uncountably many packings of discs of size $1$ and $r$, none has density larger than the claimed upper bound.
A first milestone was the proof in \cite{Flo60} that $\delta(r)$ is less than the density inside a triangle with mutually tangent discs of size $1$, $r$ and $r$ centered on its vertices,see Fig.~\ref{fig:density}, bottom center.
This upper bound was later on enhanced for $r\geq 0.6735$ in \cite{Bli69}, by proving that $\delta(r)$ is less than the density inside the union of a regular heptagon and a regular pentagon, with the heptagon (resp. pentagon) being circumscribed to a large disc (resp. small disc), see Fig.~\ref{fig:density}, bottom right (the above mentioned constant $r_B$ is the value of $r$ for which this bound reaches $\delta_1$).
In Section~\ref{sec:upper}, we explain how we obtain upper bounds by enhancing the computer-aided method used in \cite{BF22} to prove the maximal density of the $9$ triangulated packings.
In a nutshell, the interval $(0,1)$ of the possible disc size ratio is divided into sufficiently many small intervals on which a program (the \verb-C++- code is provided in the supplementary materials) searches by dichotomy an upper bound on the maximum density.
We have tried to make this paper readable as independently of \cite{BF22} as possible by recalling the main lines of the method used in \cite{BF22} (Subsec.~\ref{sec:rappel}).
It nevertheless relies on many technical details of \cite{BF22} that are omitted here and is therefore not self-contained.
The obtained upper bounds improve the previous ones for any $r\in[0.11,r_B)$ but the $8$ values in this interval which allow a triangulated binary packing.
They are depicted in Fig.~\ref{fig:density} (red curve).
As a corollary, we get:

\begin{corollary}
\label{cor:upper_bound}
One has $\delta(r)=\delta_1$ for any $r$ in
$$
\underbrace{[0.4445, 0.4532]}_{I_1} \cup
\underbrace{[0.4917, 0.5145]}_{I_2} \cup
\underbrace{[0.5666, 0.6270]}_{I_3} \cup
\underbrace{[0.6468,1)}_{I_4},
$$
where the endpoints of the $I_i$'s are exact numerical values obtained by computer.
\end{corollary}

Combining the intervals given in Cor.~\ref{cor:lower_bound} and \ref{cor:upper_bound} answers whether a given ratio $r$ allows a binary packing more dense than the hexagonal compact packing of equal discs in more than $93\%$ of the cases.
We have deliberately called both these results ``corollaries'' to emphasize the fact that the main result of this paper are the lower and upper bounds depicted in Fig.~\ref{fig:density}, which unfortunately cannot be stated in the form of a classical theorem.
The feasibility of closing the gap between the lower and upper bounds, is discussed in Section~\ref{sec:discussion}.

\section{Lower bounds by flipping and flowing}
\label{sec:lower}

\subsection{Principle}

To date, all disc packings that have been proven to maximize density (the $9$ binary ones and the hexagonal compact packing of equal discs) are found to have a triangulated contact graph, i.e. maximum number of contacts between discs \cite{FT43,BF22,Fer19,Bli69}.
This backs up the rule of thumb that the more contact between discs in a packing, the denser it is.
However, a triangulated contact graph is only possible for some very particular sizes of discs.
This is where {\em flipping and flowing} comes into play.
The principle is, starting from a particularly dense packing, to continuously modify the ratio of disc sizes while trying to keep as much contacts between discs as possible, in the hope of decreasing the density as little as possible.
The term ``flowing'' refers to the continuous modification of disc sizes.
The term ``flipping'' refers to the particular (but frequent) case where the transformation connects two packings whose graphs are triangulated and differ by one or more {\em flips}, i.e. a diagonal of a quadrilateral is replaced by the other diagonal.

\subsection{The flows}

\begin{figure}[hbt]
\centering
\includegraphics[scale=0.047]{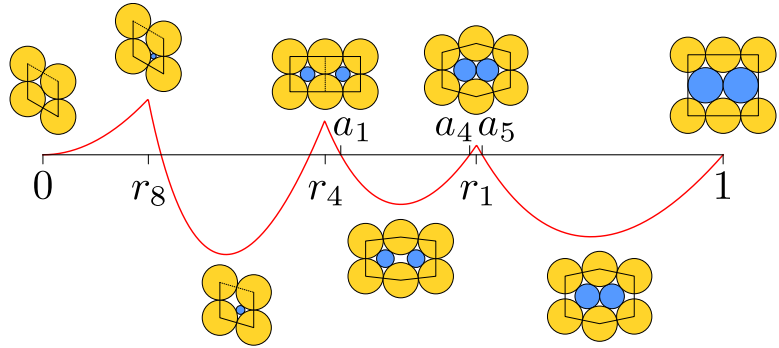}
\caption{
The flow $0\leftrightarrow r_8\leftrightarrow r_4 \leftrightarrow r_1 \leftrightarrow 1$.
}
\label{fig:flow841}
\end{figure}

\begin{figure}[hbt]
\centering
\includegraphics[scale=0.047]{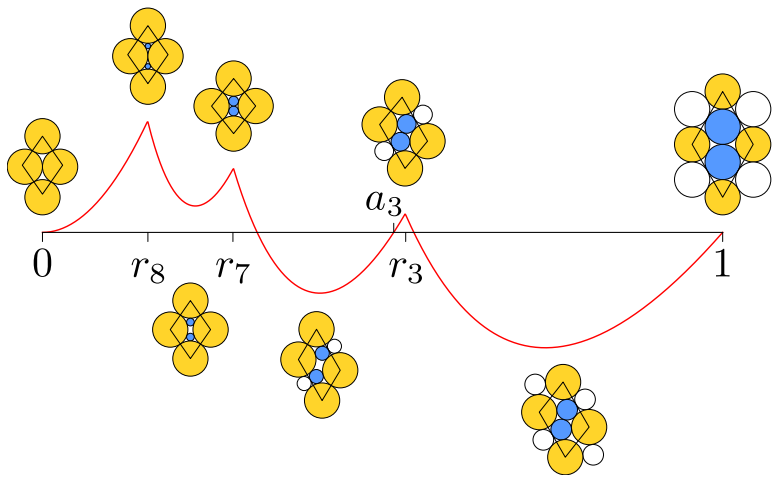}
\caption{The flow $0\leftrightarrow r_8\leftrightarrow r_7 \leftrightarrow r_3 \leftrightarrow 1$.
The white discs, not in the fundamental domain, are depicted to increase clarity.
}
\label{fig:flow873}
\end{figure}

\begin{figure}[hbt]
\centering
\includegraphics[scale=0.047]{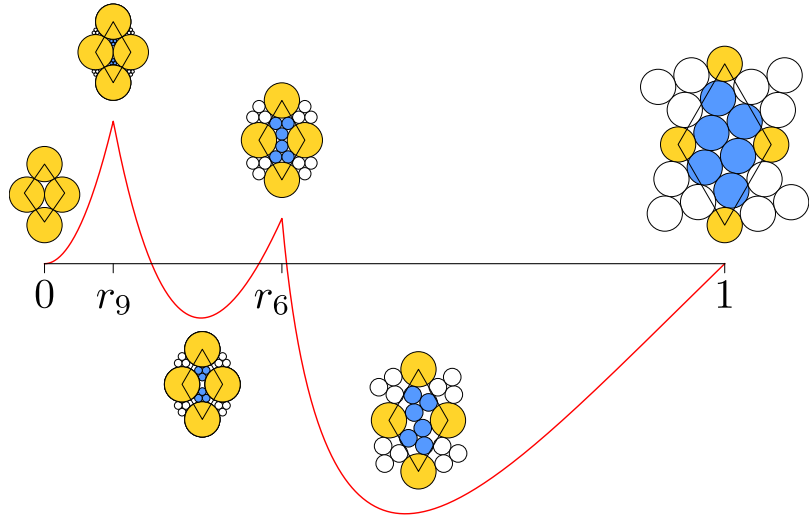}
\caption{The flow $0\leftrightarrow r_9\leftrightarrow r_6 \leftrightarrow 1$.}
\label{fig:flow96}
\end{figure}

\begin{figure}[hbt]
\centering
\includegraphics[scale=0.047]{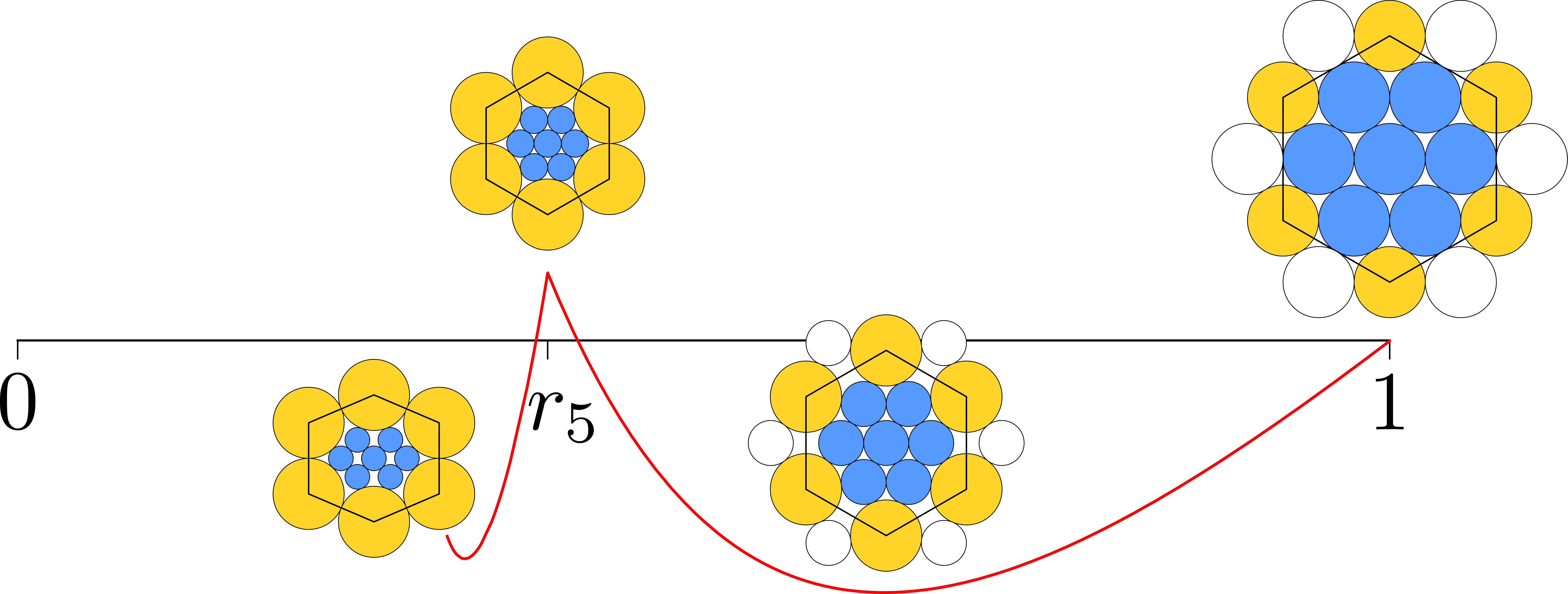}
\caption{
Flowing around $r_5$.
Among the numerous way to flow for $r<r_5$, this is the one which seems to decrease the least the density.
}
\label{fig:flow5}
\end{figure}

\begin{figure}[hbt]
\centering
\includegraphics[scale=0.047]{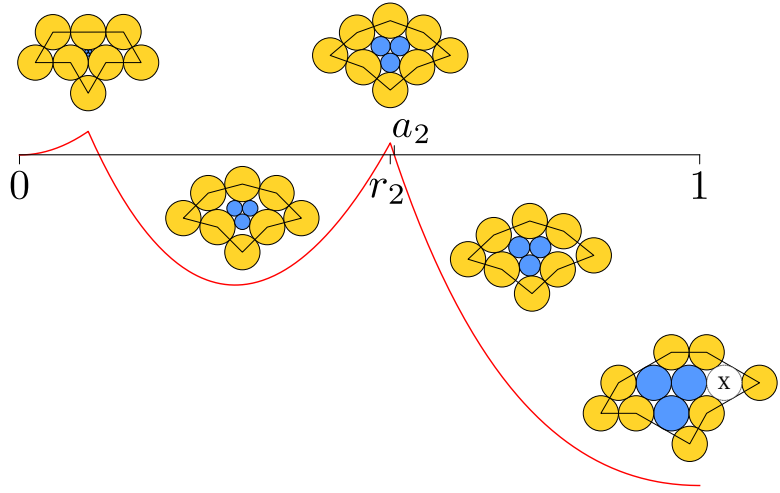}
\caption{
Flowing around $r_2$ (the fundamental domain tiles like scales).
The flow can be extended to $r=1$ but we get only $6$ out of $7$ discs of the hexagonal compact packing (the missing one is marked by an X).
}
\label{fig:flow2}
\end{figure}

\begin{figure}[hbt]
\centering
\includegraphics[scale=0.047]{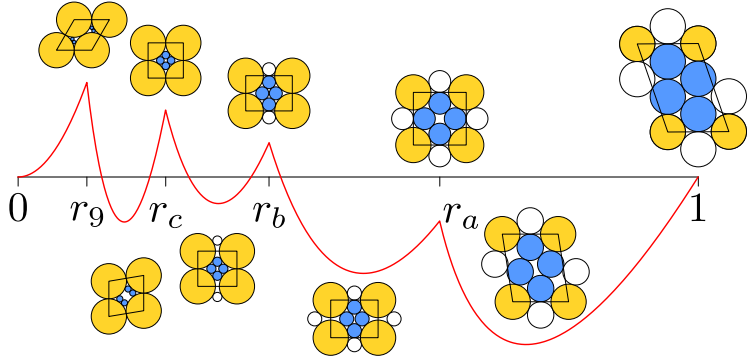}
\caption{The flow $r_9\leftrightarrow r_c \leftrightarrow r_b \leftrightarrow r_a \leftrightarrow 1$.}
\label{fig:flow9cba}
\end{figure}

Figures \ref{fig:flow841}--\ref{fig:flow9cba} describe flows between (or from) especially dense packings.
The ratios $r_i$'s and $a_i$'s are algebraic numbers whose exact values can be found in the file \verb+lower_bounds.sage+.
These packings are all triangulated, except those for $r\in\{r_a,r_b,r_c\}$ in Fig.~\ref{fig:flow9cba}.
They are all periodic and described by their fundamental domains.
The density is depicted (red curve), with the horizontal axe corresponding to the density $\delta_1$ of the hexagonal compact packing.
This yields the lower bound depicted in Fig.~\ref{fig:density}.

\begin{figure}[hbt]
\centering
\includegraphics[width=0.8\textwidth]{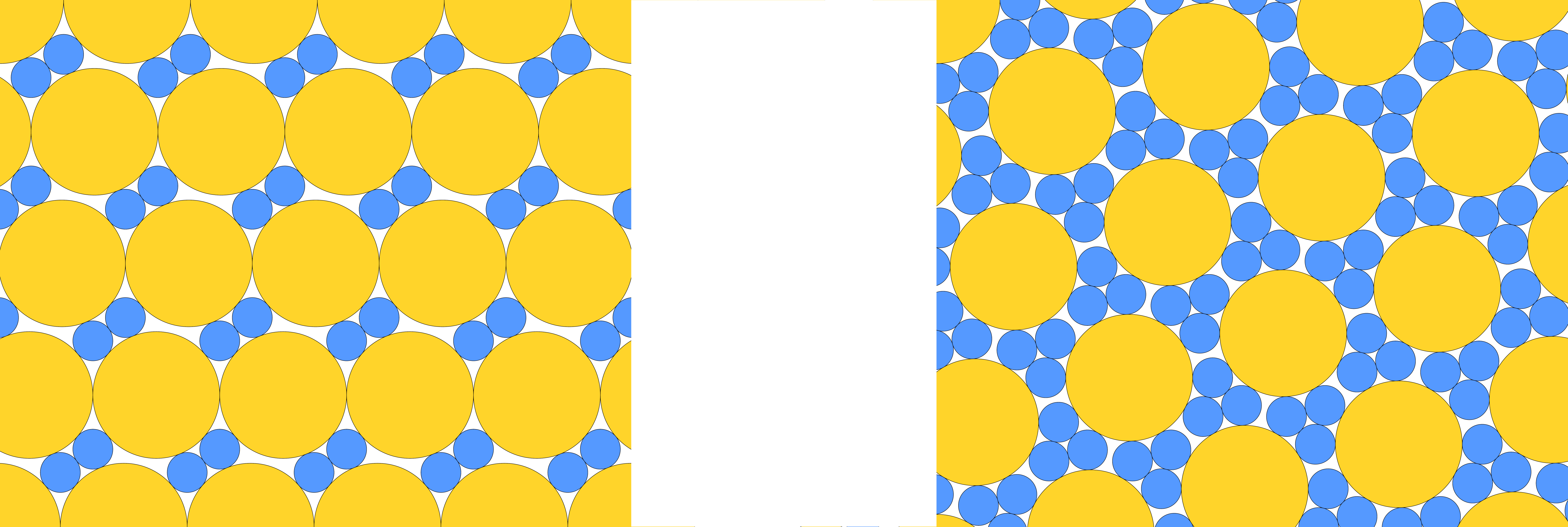}
\caption{
Two packings with the same density for $r\approx 0.3154$.
}
\label{fig:between_r6_r7}
\end{figure}

\begin{remark}
One checks that when the flow on the right of $r_7$ and the one on the left of $r_6$ cross, for $r\approx 0.3154$, the density is higher about $0.0076\%$ than $\tfrac{\pi}{2\sqrt{3}}$ (it is important for Cor.~\ref{cor:lower_bound}).
The bottom of the green valley between $r_7$ and $r_6$ in Fig.~\ref{fig:density} is thus (very slightly) above the horizontal dotted line.
For this $r$, the two flows yield packings with the same density but different structures (Fig.~\ref{fig:between_r6_r7}).
\end{remark}

\subsection{Interstitials packings}

For $r\leq r_8$, a small disc can fit in the holes of a hexagonal compact packing of large discs.
To obtain the lower bound depicted in Fig.~\ref{fig:density} for $r\leq r_8$, we simply put as much as possible discs of a hexagonal compact packing of small discs inside the holes of a hexagonal compact packing of large discs (the center of each hole is the center of a small disc).
This is not optimal (some improvements can be found in \cite{UST04}).
At least, this yields $\delta(r)>\delta_1$ over $(0,r_8)$.

\subsection{Computation}

To compute the density along a flow, we use the fact that the contacts between discs in the packings yield quadratic equations in the coordinates of disc centers and the ratio $r$.
This form a polynomial system which has to be of dimension at least $1$ to allow $r$ to vary and at most $1$ to have as much contacts as possible.
Solving this system allows to compute the positions of the discs and the density as functions of the ratio $r$.
The complete calculations are provided in the file \verb+lower_bound.sage+.
The main tool is the simple function \verb+stick((x1,y1,r1),(x2,y2,r2),r3)+, which takes the coordinates \verb+x1,y1+ of a disc $D_1$ of radius \verb+r1+, the coordinates \verb+x2,y2+ of a disc $D_2$ of radius \verb+r2+, a real number \verb+r3+ and returns the coordinates \verb+x3,y3+ (if they exist) of the disc $D_3$ of radius \verb+r3+ which is exteriorly tangent to $D_1$ and $D_2$ such that $D_1$ sees $D_3$ on the left of $D_2$.

\begin{figure}[hbt]
\centering
\includegraphics[width=0.7\textwidth]{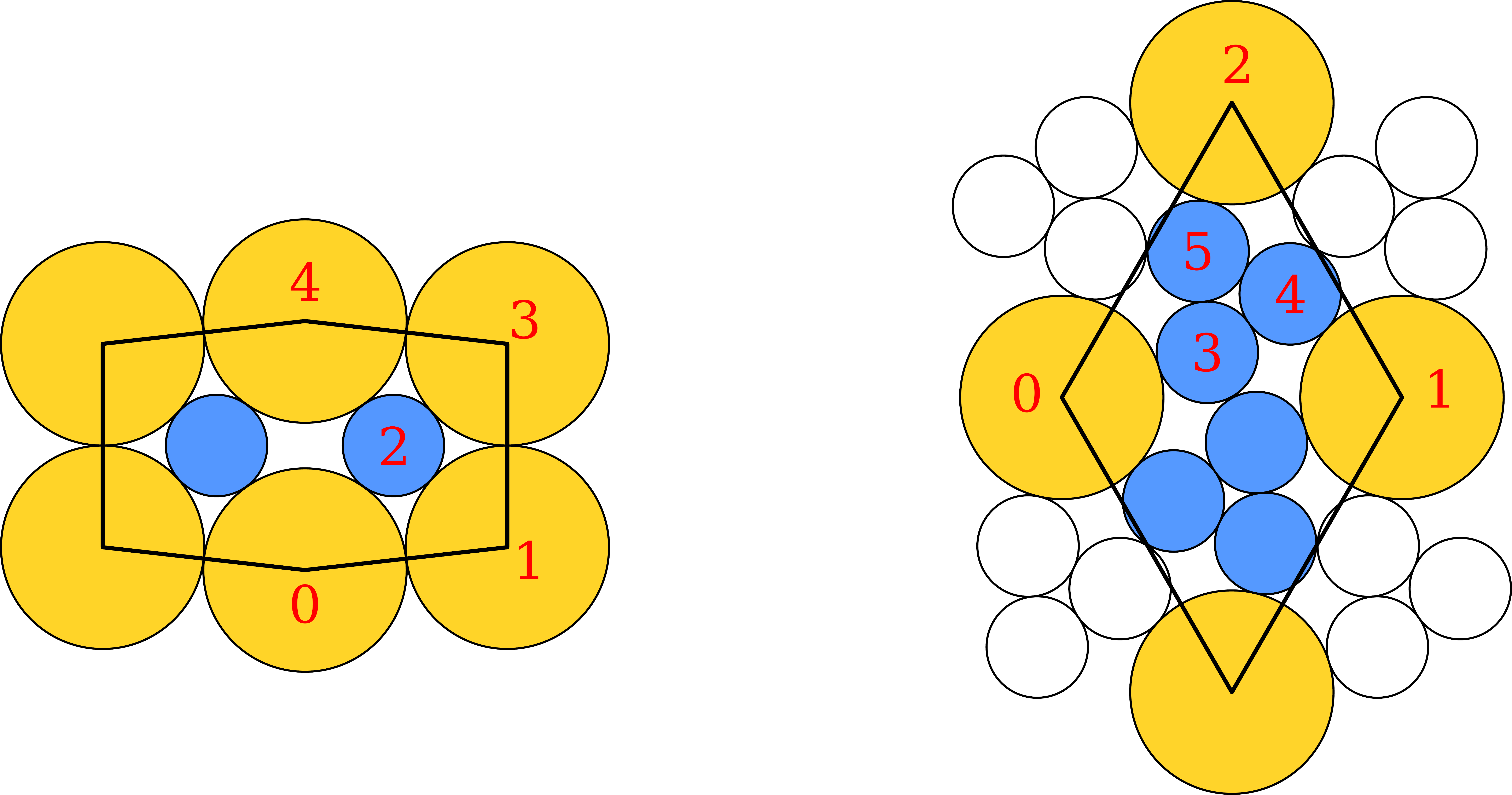}
\caption{
Disc numbering to analyze flows.
}
\label{fig:flow_ex}
\end{figure}

Consider, for example, the periodic binary packing whose fundamental domain is depicted on Fig.~\ref{fig:flow_ex}, left.
Discs are numbered and we arbitrarily set the positions of the two first discs such that they are tangent:
\begin{verbatim}
d0=(0,0,1)
d1=(2,0,1)
\end{verbatim}
We then stick one by one the three following discs:
\begin{verbatim}
d2=stick(d0,d1,r)
d3=stick(d2,d1,1)
d4=stick(d2,d3,1)
\end{verbatim}
This allows to compute the density in the fundamental domain:
\begin{verbatim}
d=(2+2*r^2)/((d4[1]+(d3[1]-d1[1]))*d1[0])
\end{verbatim}
This yields for the density along the flow between $r_4$ and $r_1$ the expression
$$
\frac{\pi(r^2 + 1)(r + 1)^4}{16(r + 2)\sqrt{r+2}r\sqrt{r}}.
$$
All the considered cases are similar, except two which are slightly more complicated.
The first case is the periodic binary packing whose fundamental domain is depicted on Fig.~\ref{fig:flow_ex}, right.
It corresponds,  in the previous subsection, to the flow between $r_6$ and $1$.
In this case, it is not possible to describe the packing by sticking discs one by one to two previous discs.
To get around this problem, we define a multivariate polynomial ring over $\mathbb{Q}$ whose variables are the coordinates of the disc centers, the disc ratio \verb+r+ and the density \verb+d+:
\begin{verbatim}
K.<x1,y1,x2,y2,x3,y3,x4,y4,x5,y5,r,d>=QQ[]
\end{verbatim}
We then add the equations which describe the disc contacts or symmetries of the packing (each polynomial must be equal to zero):
\begin{verbatim}
eqs=[y1,
x1-2*x2, 
x1+x2-(x3+x4+x5), 
y1+y2-(y3+y4+y5), 
x3^2+y3^2-(1+r)^2, 
(x4-x3)^2+(y4-y3)^2-(r+r)^2, 
(x4-x5)^2+(y4-y5)^2-(r+r)^2, 
(x3-x5)^2+(y3-y5)^2-(r+r)^2, 
(x4-x1)^2+(y4-y1)^2-(1+r)^2, 
(x2-x5)^2+(y2-y5)^2-(1+r)^2, 
(x1-2*x3)^2+(y3+y3)^2-(r+r)^2,
d*x1*y2-(1+6*r^2)]
\end{verbatim}
We then compute the ideal defined by these equations and eliminate all the variables but \verb+r+ and \verb+d+:
\begin{verbatim}
I=ideal(K,eqs)
J=I.elimination_ideal([x1,y1,x2,y2,x3,y3,x4,y4,x5,y5])
P=QQ[r,d](J.gens()[0])
\end{verbatim}
This yields a polynomial of degree $9$ in \verb+r+ and $6$ in \verb+d+.
Actually, we can remove a factor which would correspond to a density larger than $1$.
We get an irreducible polynomial of degree $2$ in \verb+d^2+ whose coefficients are polynomials in \verb+r+.
This allows to get a closed-form expression for the density along the flow between $r_6$ and $1$:
$$
\frac{\pi(6  r^2 + 1)\sqrt{47  r^4 + 84  r^3  + 54  r^2 + 12  r + 3 \textrm{$-$} (7  r^3 + 13 r^2 + 9 r + 3 ) \sqrt{45r^2 \textrm{$-$} 6r \textrm{$-$} 3}}}{\sqrt{6}(r^4 + 12  r^2 + 12  r + 3)}.
$$
The second case corresponds, in the previous subsection, to the flow on the left of $r_5$.
It is similar, except that we eventually get a polynomial in degree $6$ in \verb+d^2+: we cannot derive a closed-form expression, but an implicit plot is possible.

\begin{remark}
The ``algebraic'' method by ideal elimination can actually be used for all the flows.
The more ``pedestrian'' method of adding the discs one by one leads however to much faster calculations in practice.
\end{remark}

\section{Upper bounds via localizing potentials}
\label{sec:upper}

\subsection{Checking an upper bound for a fixed ratio}
\label{sec:rappel}

Given a ratio $r$ and a candidate upper bound $\delta$ on the maximum density, we want to check whether $\delta(r)\leq\delta$ or not.
The strategy used is the one in \cite{BF22}, which resembles the one used by Hales to prove the Kepler conjecture, nicely exposed in \cite{Lag02}.
Here we will only sketch this strategy and omit all technical details that are not strictly necessary to understand what follows.
We refer to \cite{BF22} for the complete description (which is unfortunately necessary to fully understand the \verb-C++- code provided in the supplementary materials).

Given a packing, consider a triangulation of its disc center (we used the Fejes-Moln\'ar triangulation \cite{FM58} for its suitable properties).
The density inside the triangles, that is, the proportion of the triangle area covered by the discs of the packing, can greatly vary from one triangle to the other.
In particular, it may be larger than the candidate upper bound $\delta$ on the maximum density for any packing.
The idea is to prove that whenever a triangle is too dense, there are necessarily neighbor triangles whose densities are low enough so that the density drops below $\delta$ once averaged over all these triangles.
More precisely, the density of each triangle is distributed among its three vertices, and it is proven that no matter how a vertex of the triangulation of a packing is surrounded, the average of the densities distributed to this vertex by the triangles to which it belongs is at most $\delta$.
Key parameters appear to be the so-called ``base vertex potentials'' $V_{111}$, $V_{11r}$, $V_{r1r}$, $V_{1r1}$, $V_{1rr}$ and $V_{rrr}$.
Each $V_{ijk}$ specifies, in a triangle whose vertices are center of mutually tangent discs of radii $i$, $j$ and $k$, the amount of density distributed to the center of the disc of radius $j$.
Since the density of such a triangle is determined and since there are $4$ different such triangles (depending on the radii $111$, $11r$, $1rr$ or $rrr$ of the discs centered on vertices), two of the $6$ base vertex potentials have to be chosen and the others follow.
This choice has a strong influence on whether the maximum density can be proven to be less than $\delta$.
This choice turned out to be relatively simple to make in each of the 9 ratios considered in \cite{BF22}.
But here the ratio can be any, so it is necessary to automate this choice in a satisfactory way.

To explain how this choice is made, we need to mention two other parameters that appear in \cite{BF22}, namely $m_1$ and $m_r$.
The base vertex potential $V_{ijk}$ is indeed defined only in triangles with mutually tangent discs.
For other triangles, $V_{ijk}$ grows linearly with the angle in $j$ with a proportionality factor $m_j$.
When the angle in $j$ increases, the triangle tends to be less dense and can ``absorb'' more density from its neighbors, which encourages taking large values of $m_j$ and $V_{ijk}$.
However, the total density that each triangle can absorb cannot is limited, which in turns limits the value that $m_j$ and $V_{ijk}$ can take.
This can be formalized by a list of inequalities on $m_1$, $m_r$ and two of the $V_{ijk}$'s (since the $4$ other ones follow as already mentioned).
In the program \verb+upper_bound.cpp+, we consider $V_{11r}$ and $V_{1rr}$.
The inequalies are defined in the function \verb+add_vertex_positivity_constraints+ of the program \verb+parameters_xy.cpp+.
This yields a $4$-dimensional {\em parameter polytope} (polyhedron \verb+P+ in the function \verb+find_delta+ in \verb+upper_bound.cpp+) in which a point has to be chosen.

The following automatic choice was developed through trial and error, with the goal being to prove an upper bound $\delta$ as low as possible.
It is implemented in the function \verb+set_xy_generic+ in the program \verb+upper_bound.cpp+.
First, it sets the ratio $m_1/m_r$ to be equal to the maximal possible value of $m_1$ divided by the maximal possible value of $m_r$.
This amounts to intersecting the parameter polytope with a hyperplane to get a new polytope.
We then take either the barycenter of the vertices of this new polytope if $r>0.55$, or a vertex of it which minimizes $m_1$ otherwise.
This defines the parameters $V_{11r}$ and $V_{1rr}$.
We do not consider the values $m_1$ and $m_r$ it defines because our program deals with rational polytopes whereas $r$ can be not rational (actually, it can even be an interval as we shall later see).
Instead, we keep only the numerical values of $V_{11r}$ and $V_{1rr}$ and then proceed as in \cite{BF22} to compute $m_1$ and $m_r$ (as well as the parameter $\varepsilon$, not mentioned here) which will decide exactly whether $\delta$ is a suitable upper bound or not.

\subsection{Finding an upper bound for an interval of ratios}

The previous subsection explained how to check a given upper bound for a given ratio.
But the goal of this paper is to find an upper bound as good as possible for any ratio.
This yields two problems:
\begin{enumerate}
\item we have no candidate upper bound on the maximal density;
\item we have a continuum of ratios to consider.
\end{enumerate}

Assuming $r$ is fixed, the first problem can be fixed by dichotomy on the candidate density.
Namely, we maintain two variables $\delta_\textrm{inf}$ and $\delta_\textrm{sup}$, such that the proof that we have an upper bound on the exact maximal density succeeds for $\delta_\textrm{sup}$ but fails for $\delta_\textrm{inf}$.
We start with $\delta_\textrm{inf}$ slightly less than $\delta_1$ and $\delta_\textrm{sup}$ equal to the Florian upper bound.
At each step, we check whether the proof works for $\tfrac{1}{2}(\delta_\textrm{inf}+\delta_\textrm{sup})$ and update $\delta_\textrm{inf}$ and $\delta_\textrm{sup}$ accordingly.
We stop when $\delta_\textrm{sup}-\delta_\textrm{inf}$ is smaller than a fixed precision, namely $0.0001$.
We finally output $\delta_\textrm{sup}$: it is an upper bound on the maximal density (and the best that we get by our method, up to the fixed precision).
This is done in the function \verb+find_delta+ in \verb+upper_bound.cpp+.

To fix the second problem, we can subdivide $(0,1)$ is many small intervals and rely on interval arithmetic to consider each small interval as a value for $r$.
The smaller the intervals, the better the precision, thus the better the upper bound on the maximal density.
We thus want to subdivide $(0,1)$ into intervals as small as possible, with the limiting factor being the computation time.

Actually, we found it more relevant to compute upper bound only for regularly spaced discrete values of $r$ with maximal precision.
This indeed gives a fair idea of which bounds could be obtained because the maximal density is quite regular:

\begin{proposition}
\label{prop:regularity}
For $x<y$ in $[0,1]$, the maximal density satisfies
$$
\frac{|\delta(y)-\delta(x)|}{|y-x|} \leq \frac{\pi}{y^2\sqrt{3}}.
$$
\end{proposition}

\begin{figure}[hbt]
\centering
\includegraphics[width=1\textwidth]{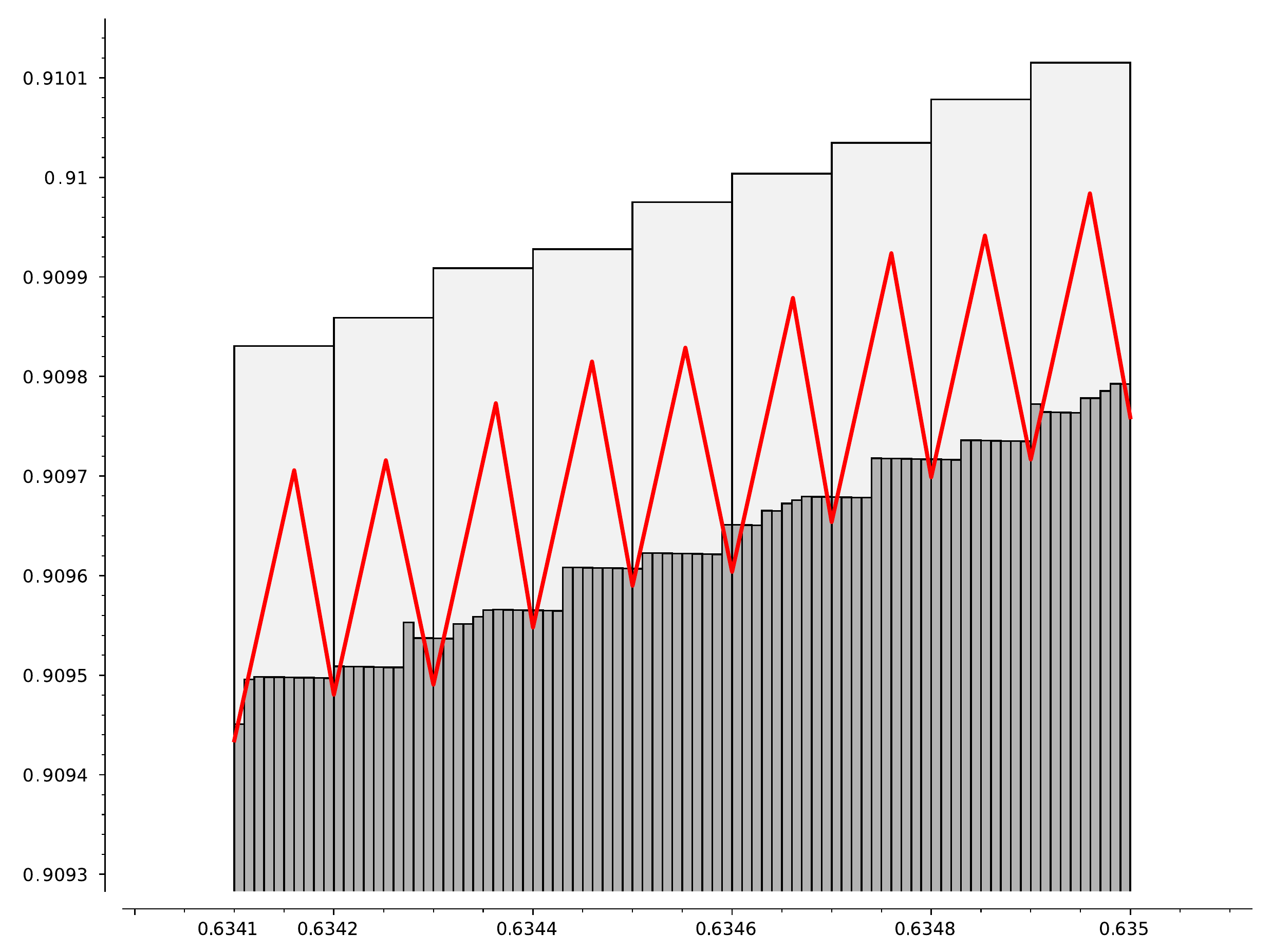}
\caption{
The upper bound obtained with subintervals of length $0.0001$ (resp. $0.00001$) in light gray (resp. dark gray).
The local minima of the red broken line are the upper bounds obtained with the corresponding value of $r$, and the line itself is deduced from Prop.~\ref{prop:regularity}.
Since connecting the minima by segments is expected to give a fair idea of the best upper bound that can be obtained by our method, we only represented these minima in Fig.~\ref{fig:density}, \ref{fig:J34}, \ref{fig:J23}, \ref{fig:J12} and \ref{fig:J0}.
}
\label{fig:precision}
\end{figure}

The proof of this proposition is given in Appendix~\ref{app:regularity}.
In particular, the maximal density is $\tfrac{\pi}{a^2\sqrt{3}}$-Lipshitz over any interval $[a,b]\subset(0,1]$.
This makes the computation much lighter.
For example, on our computer, subdividing the interval $[0.6269, 0.6469]$ takes $27 $min. for subintervals of length $0.0001$ and $1930$ min. ($32$h.) for subintervals of length $0.00001$.
In comparison, discrete values of $r$ with a step of $0.0001$ takes $42$ min. and yields for these discrete points an upper bound on the maximal density which is comparable to the one obtained with intervals of length $0.00001$, in around $50$ times less computation time.
The regularity allows to extend this upper bound everywhere.
Fig.~\ref{fig:precision} illustrates this point.

Table~\ref{tab:discrete_bounds} gives the step and execution time on each connected component of the complement of $I_1\cup I_2\cup I_3\cup I_4$ (the intervals on which the maximal density has been proven to be $\delta_1$).
A list of values can be found in the supplementary materials (file \verb+trace_upper_bound.txt+), see also Fig.~\ref{fig:J34}, \ref{fig:J23} and \ref{fig:J12} for a zoom of Fig.~\ref{fig:density}.

\begin{table}[hbt]
\centering
\begin{tabular}{|c|c|c|}
\hline
Interval & step & execution time\\
\hline
$[0.6269, 0.6469]$ & $0.0001$ & 42 min.\\
$[0.5145, 0.5666]$ & $0.0001$ & 98 min.\\
$[0.4532, 0.4917]$ & $0.0001$ & 97 min.\\ 
$[0.106, 0.445]$ & $0.001$ & 16h.\\
\hline
\end{tabular}
\caption{
Interval, step between consecutive values of $r$ and execution time.
}
\label{tab:discrete_bounds}
\end{table}

\subsection{Checking $\delta(r)=\delta_1$ over $I_1$, $I_2$, $I_3$ and $I_4$}

The upper bound obtained in the previous subsection suggest $\delta(r)=\delta_1$ over the intervals $I_1$, $I_2$, $I_3$ and $I_4$ defined in Corollary~\ref{cor:upper_bound}.
To prove this, it is no longer enough to consider discrete values of r: we really have to use interval arithmetic.
Instead of partitionning each $I_i$ in many equally small intervals, we use again the dichotomy to subdivide as little as possible.
Namely, we bisect each $I_i$ while the precision does not suffice to conclude.
But we no longer need to make a dichotomy on the candidate density since it is $\delta_1$.
Table~\ref{tab:subintervals} gives the total number of subintervals into which the program \verb+upper_bound_HCP.cpp+ eventually divided each interval $I_i$ as well as the execution time on our Laptop (i5-7300U CPU 2.60GHz).

\begin{table}[hbt]
\centering
\begin{tabular}{|c|c|c|}
\hline
Interval & subintervals & execution time\\
\hline
$I_1$ & $259$ & 30 min.\\
$I_2$ & $318$ & 20 min.\\
$I_3$ & $379$ & 32 min.\\
$I_4$ & $845$ & 20 min.\\
\hline
\end{tabular}
\caption{
Interval, number of subintervals and execution time.
}
\label{tab:subintervals}
\end{table}

\section{Discussion}
\label{sec:discussion}

\begin{figure}[hbt]
\centering
\includegraphics[width=\textwidth]{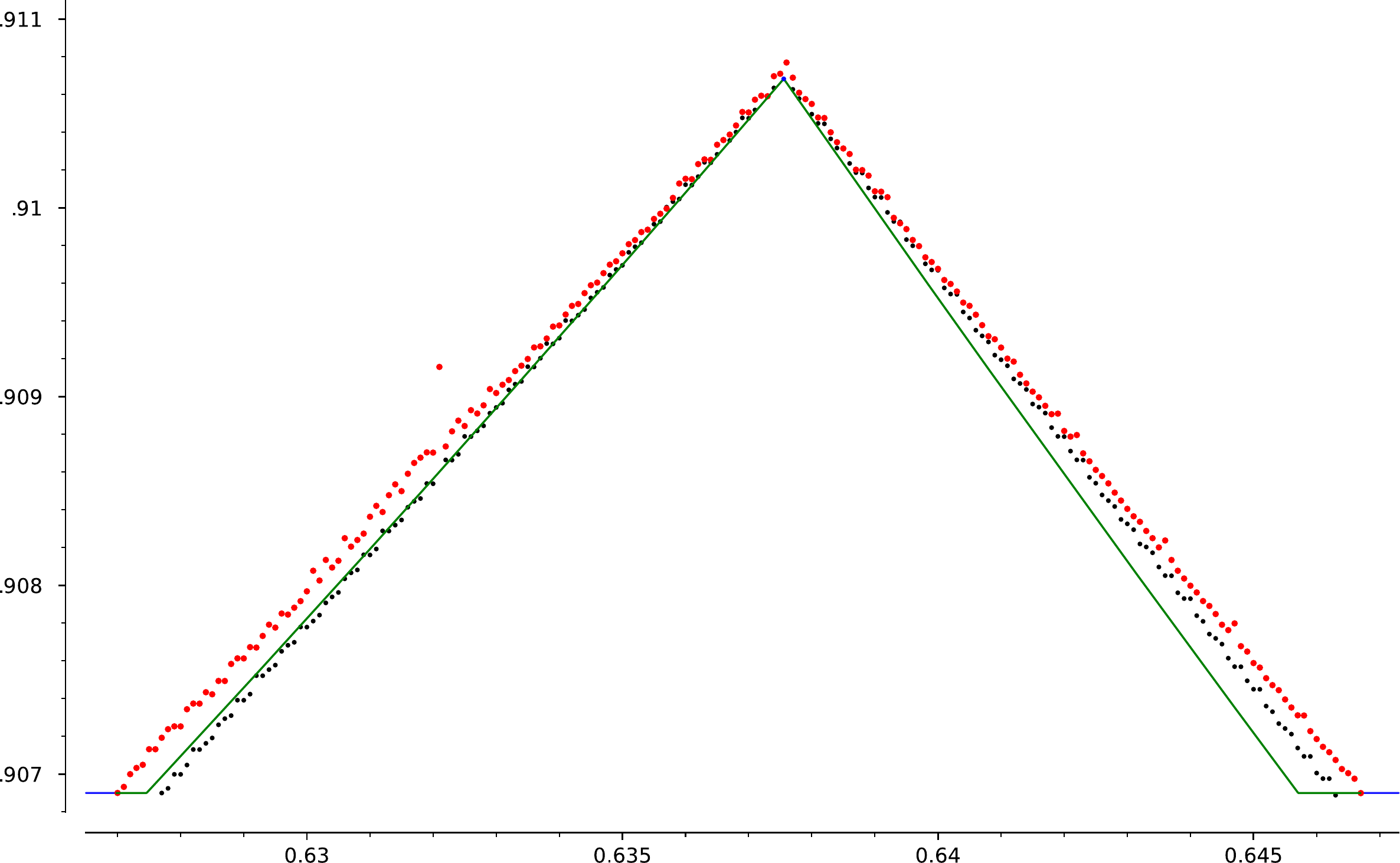}
\caption{
Between $I_3$ and $I_4$, that is, for $r\in[0.6269, 0.6469]$.
The upper bound (red points) is quite close to the lower bound (green curve).
We conjecture that the lower bound is tight on this interval.
For $r_1\approx 0.6375$, the lower bound reaches its maximum and has been proven to be tight \cite{BF22}.
}
\label{fig:J34}
\end{figure}

\begin{figure}[hbt]
\centering
\includegraphics[width=\textwidth]{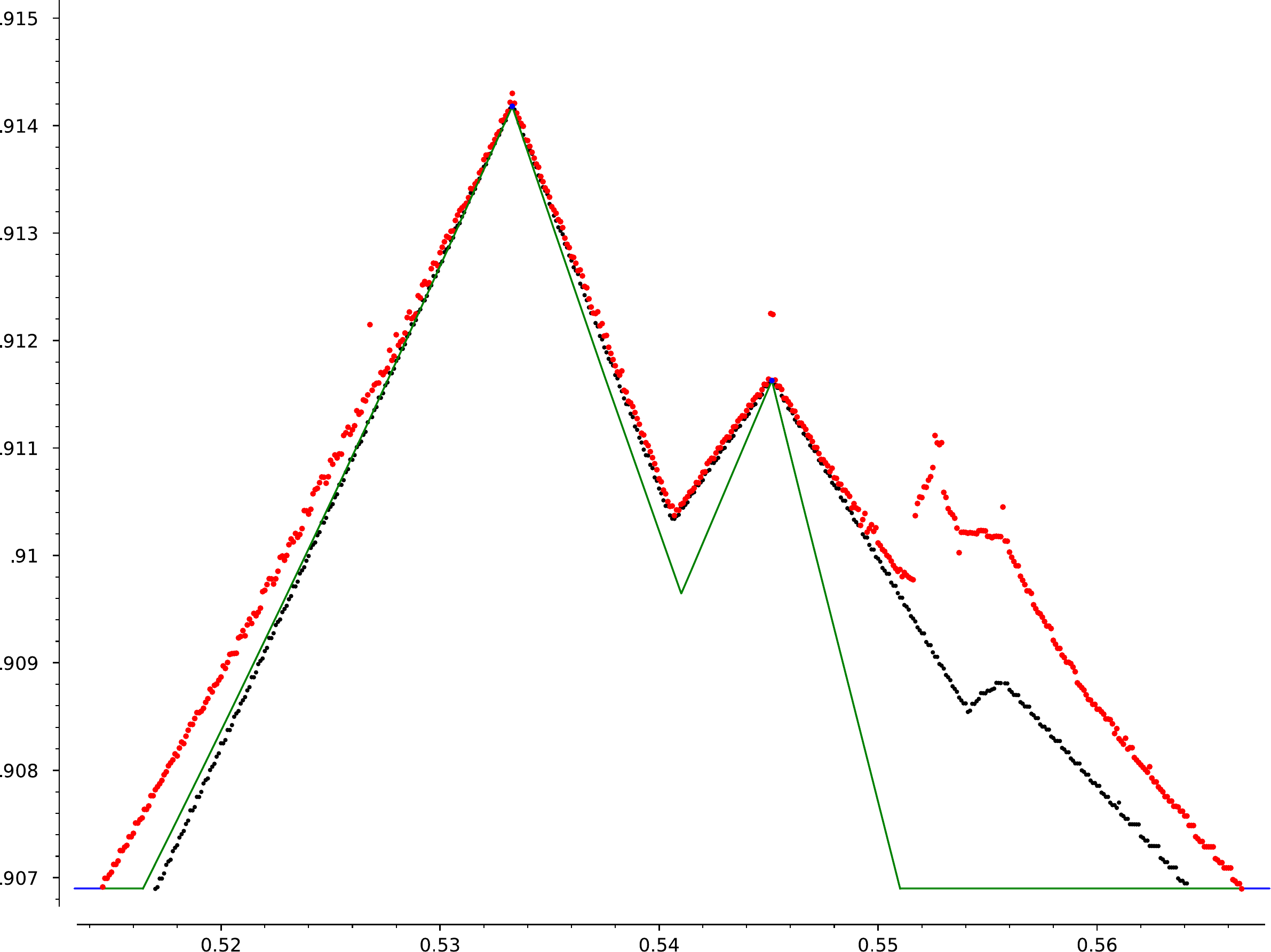}
\caption{
Between $I_2$ and $I_3$, that is, for $r\in [0.5145, 0.5666]$.
For $r_3\approx 0.5332$ and $r_2\approx 0.5451$, the lower bound reaches two local maxima and has been proven to be tight \cite{BF22}.
For $r\leq r_2$, the upper bound (red points) are quite close to the lower bound (green curve).
We conjecture that the lower bound is tight on this interval.
For $r\geq r_2$ the lower and upper bounds diverge.
The difference between black and red dots is moreover relatively large between $0.552$ and $0.558$, suggesting (but not proving) that the choice of parameters could be optimized.
We do not exclude the possibility that the lower bound is not optimal, i.e. that a better flow could be defined to the right of $r_2$ (remind Fig.~\ref{fig:flow2}).
}
\label{fig:J23}
\end{figure}

\begin{figure}[hbt]
\centering
\includegraphics[width=\textwidth]{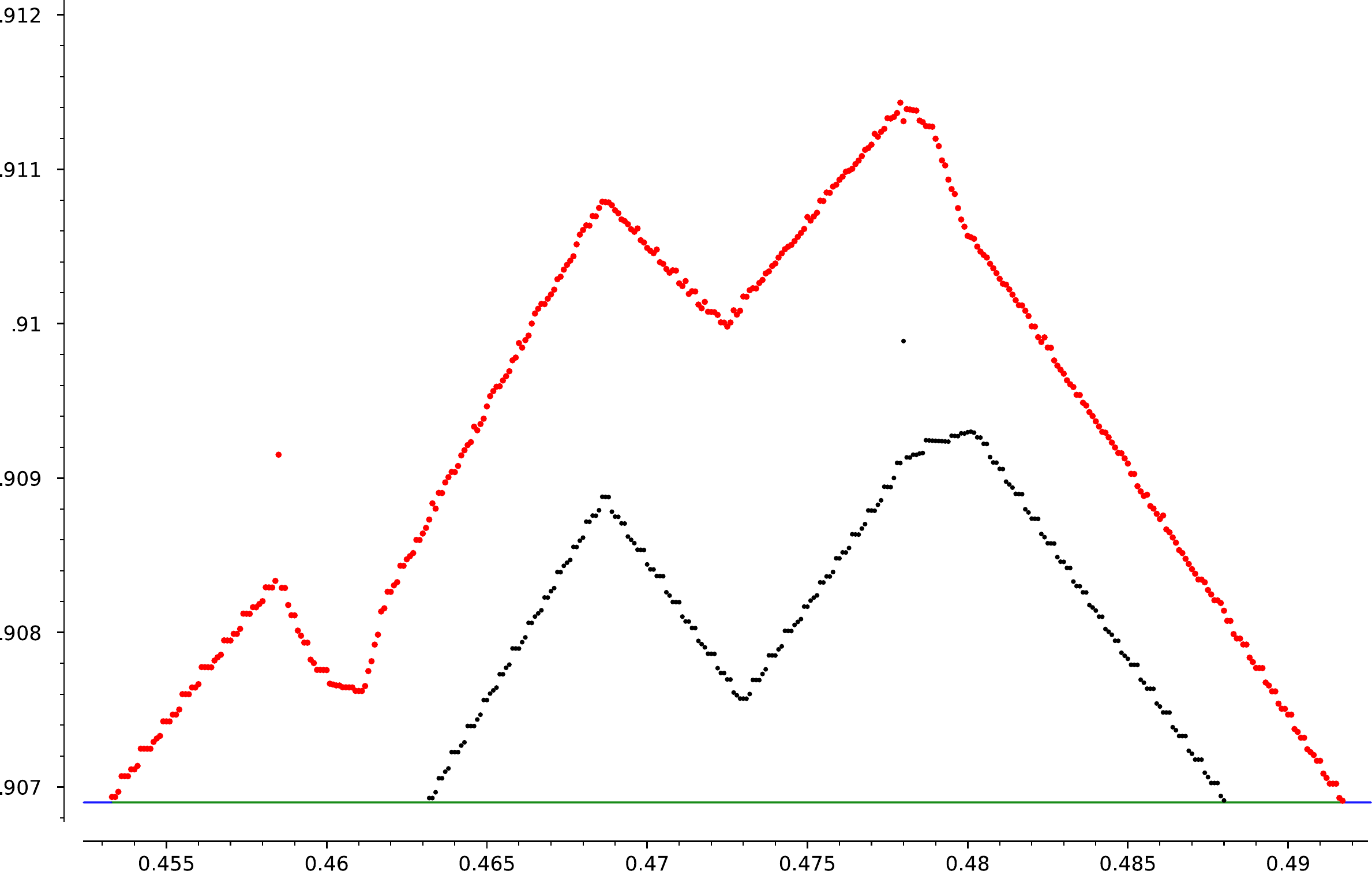}
\caption{
Between $I_1$ and $I_2$, that is, for $r\in [0.4532, 0.4917]$, there is a sort of ``mysterious island''.
The lower bound is indeed $\delta_1$ over this whole interval: no packing more dense than the hexagonal compact packing is known.
The difference between red and black points suggest that the choice of parameters could be optimized, but does not leave any hope of reducing the upper bound to $\delta_1$ over the whole interval.
Our conjecture is that the lower bound is tight, that is, this mysterious island is an artefact of our method.
This possible ``artefact'' is discussed in more detail in the text.
}
\label{fig:J12}
\end{figure}

\begin{figure}[hbt]
\centering
\includegraphics[width=\textwidth]{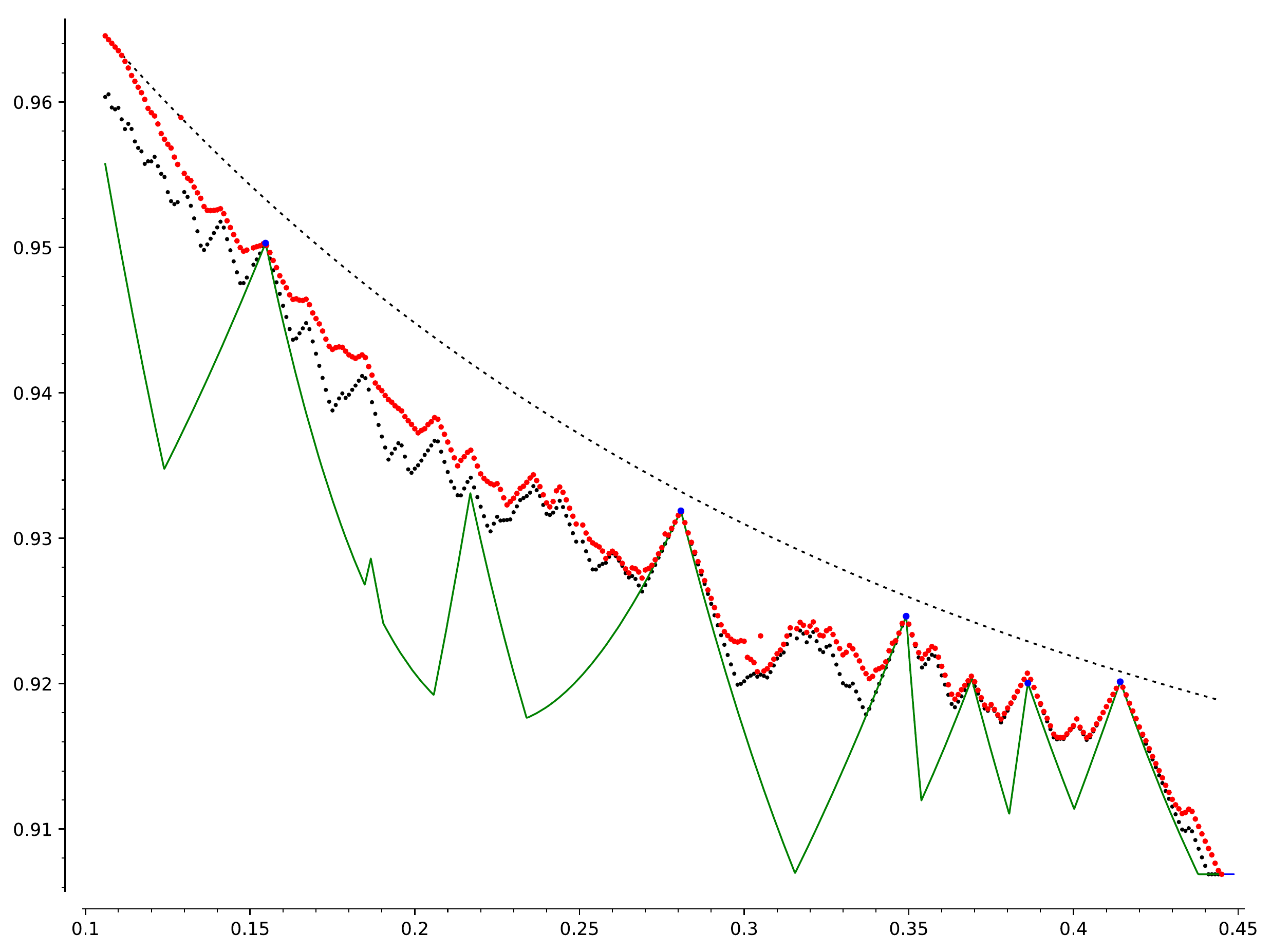}
\caption{
On the left of $I_1$, that is, $r\in [0.106, 0.445]$.
The lower bound has been proved to be tight on its local maxima at $r_8\approx 1547$, $r_7\approx 0.2808$, $r_6\approx 0.3292$, $r_5\approx 0.3861$ and $r_4\approx 0.4142$ \cite{BF22}.
We conjecture that it still holds on a neighborhood of these ratios.
We also conjecture that the lower bound is tight at its local maximum at $r_b\approx 0.3691$.
This is less clear around $r_c\approx 0.2168$.
Upper and lower bound are quite different in the valleys betweens peaks of the lower bound.
This can be due to artefacts (this is our hypothesis - see text) or to unknown dense packings.
}
\label{fig:J0}
\end{figure}

Figures~\ref{fig:J34}, \ref{fig:J23}, \ref{fig:J12} and \ref{fig:J0} depicts lower and upper bound on the maximal density over each interval of interest.
The upper bounds are depicted by red points and the lower bound by a green curve.

The black points indicate the smallest $\delta$ that we get if we only check that the parameter polytope is not empty (this corresponds to the first pass in the function \verb+find_delta+ in \verb+upper_bound.cpp+).
These points thus give a lower bound on the best upper bound that we can hope to obtain with our method, assuming that we know how to choose the parameters in an optimal way.
This lower bound is however not necessarily achievable (this is in particular the case when the black points are under the green curve, i.e., under the lower bound proved on the maximum density): the parameter polytope could be not empty without containing any parameter that allows to prove the candidate bound. 
Indeed, the polytope is not calculated in an exact way, whereas the verification once the parameters are chosen is exact and it is used to establish the upper bound.

Nevertheless, for the ratios with a large difference between red and black points, the upper bound can probably be improved by choosing the parameters more finely.
The upper bound turns out to be very sensitive to the choice of parameters in the polytope, which is therefore a delicate exercise (especially since it is done automatically here because of the very large number of different ratios considered).

The black points also show that the largest intervals over which one can expect our method to prove that the maximal density is $\delta_1$ are only slightly larger than the $I_i$'s in Corollary~\ref{cor:upper_bound}, namely:
$$
[0.4398, 0.4644],\qquad
[0.4862,0.5182],\qquad
[0.5624,0.6285],\qquad
[0.6455,1].
$$

\begin{figure}[hbt]
\centering
\includegraphics[width=\textwidth]{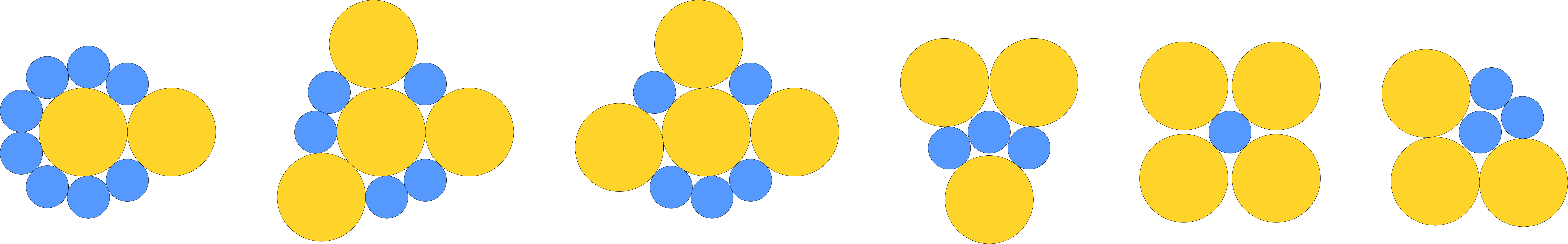}
\caption{
Discs of radius $1$ and $0.48$.
The large disc can be surrounded by discs to form a pattern which is locally quite optimal in terms of density (three leftmost patterns).
This is also the case, to a lesser extent, of the small discs (three rightmost patterns).
However, it seems impossible to combine these patterns to form a dense packing in the plane.
We can indeed start from one of these patterns, but the more we add discs, the less the local patterns can resemble those represented here.
}
\label{fig:48}
\end{figure}

In Fig.~\ref{fig:J12} and \ref{fig:J0}, we used the term ``artefact''.
By this, we mean a peak in the upper bound which is thought to be well above the exact maximal density.
Such peaks appear for ratios such that the discs fit together particularly well around one disc.
Indeed, the method developed in \cite{BF22} is rather ``local'': it distributes the densities between each disc and its close neighbors and bounds from above the resulting average density.
However, these locally dense arrangements may not combine well on a more global scale, leading to packings in the whole plane that are actually much less dense than the obtained upper bound.
Figure \ref{fig:48} illustrates the case $r=0.48$.
The smaller $r$, the more often discs fit together particularly well around one disc.
This explain why there are more and more peaks in the red or black curves in Fig.~\ref{fig:density} or \ref{fig:J0} when $r$ becomes small.
To get around this problem, it will probably be necessary to modify the method to make it less local, i.e., to distribute the densities on a larger scale. This unfortunately makes the method even more complex.

\appendix
\section{Proof of Prop.~\ref{prop:regularity}}
\label{app:regularity}

\begin{proof}
Recall that the maximal density can be approached arbitrarily close by the density of a periodic packing.
Fix $\varepsilon>0$.

First, consider a periodic packing of discs of radii $y$ and $1$ and density at least $\delta(y)-\varepsilon$.
Assume its fundamental domain has area $A$ and contains $p$ discs of radius $y$ and $q$ discs of radius $1$.
Hence
$$
\frac{p\pi y^2+q\pi}{A}\geq \delta(y)-\varepsilon.
$$
By replacing each disc of radius $y$ by a smaller disc of radius $x$ with the same center, we get a packing of discs of radii $x$ and $1$ whose density is
$$
\frac{p\pi x^2+q\pi}{A}.
$$
Since this density is, by definition, less or equal than $\delta(x)$, we have
$$
\delta(y)-\varepsilon \leq \frac{p\pi y^2+q\pi}{A}\leq\delta(x)+\frac{p\pi(y^2-x^2)}{A}=\delta(x)+\frac{p\pi(x+y)}{A}(y-x).
$$
In the initially considered packing, the fraction of $A$ covered by the discs of radius $y$ is at most $\tfrac{\pi}{2\sqrt{3}}$, the maximal density of a packing by equal discs.
Hence
$$
\frac{p\pi y^2}{A}\leq\frac{\pi}{2\sqrt{3}}.
$$
This yields
$$
\delta(y)-\varepsilon \leq \delta(x)+\frac{\pi(x+y)}{2y^2\sqrt{3}}(y-x) \leq \frac{\pi}{y^2\sqrt{3}}(y-x).
$$
Taking $\varepsilon\to 0$ gives the first half of the claimed inequality.

Conversely, consider a periodic packing of discs of radii $x$ and $1$ and density at least $\delta(x)-\varepsilon$.
Assume its fundamental domain has area $A$ and contains $p$ discs of radius $x$ and $q$ discs of radius $1$.
Hence
$$
\frac{p\pi x^2+q\pi}{A}\geq \delta(x)-\varepsilon.
$$
By scaling the whole packing by $y/x$, we get a packing of discs of radii $y$ and $y/x>1$ whose fundamental domain has area $A\times(y/x)^2$.
Then, by replacing each disc of radius $y/x$ by a smaller disc of radius $1$ with the same center, we get a packing of discs of radii $y$ and $1$ whose density is
$$
\frac{p\pi y^2+q\pi}{A(y/x)^2}=\frac{p\pi x^2}{A}+\frac{q\pi x^2}{Ay^2}.
$$
Since this density is, by definition, less or equal than $\delta(y)$, we have
$$
\delta(x)-\varepsilon
\leq \frac{p\pi x^2+q\pi}{A}
=\frac{p\pi y^2+q\pi}{A(y/x)^2}+\frac{q\pi}{A}-\frac{q\pi x^2}{Ay^2}
\leq\delta(y)+\frac{q\pi (x+y)}{Ay^2}(y-x).
$$
In the initially considered packing, the fraction of $A$ covered by the discs of radius $1$ is at most $\tfrac{\pi}{2\sqrt{3}}$, the maximal density of a packing by equal discs.
Hence
$$
\frac{q\pi}{A}\leq\frac{\pi}{2\sqrt{3}}.
$$
This yields
$$
\delta(x)-\varepsilon
\leq \delta(y)+\frac{\pi(x+y)}{2y^2\sqrt{3}}(y-x)
\leq \delta(y)+\frac{\pi}{y^2\sqrt{3}}(y-x).
$$
Taking $\varepsilon\to 0$ gives the second half of the claimed inequality.
\end{proof}

\bibliographystyle{alpha}
\bibliography{bounds}

\end{document}